\documentclass{amsart}

\usepackage[latin1]{inputenc}
\usepackage{subfigure}
\usepackage{pstricks}
\usepackage{fp}
\usepackage{multido}

\usepackage{graphicx}

\usepackage[round]{natbib}

\newtheorem{theo}{Theorem}

\newtheorem{lemm}[theo]{Lemma}

\newcommand{\pref}[1]{(\ref{#1})}

\def\cf{{\it cf. }}
\def\ie{{\it i.e. }}

\def\si{\sigma}

\def\bx{\bar x}

\def\ba{\bar a}

\newcommand{\flecheD}{\psline{->}(0,0)(.75,0)\psline(0,0)(1,0)}
\newcommand{\flecheG}{\psline{->}(0,0)(-.75,0)\psline(0,0)(-1,0)}
\newcommand{\flecheH}{\psline{->}(0,0)(0,.75)\psline(0,0)(0,1)}
\newcommand{\flecheB}{\psline{->}(0,0)(0,-.75)\psline(0,0)(0,-1)}

\newcounter{ISL}
\newcounter{ISH}
\newcommand{\IceGrid}[2]{%
\setcounter{ISL}{#1}\addtocounter{ISL}{-1}%
\setcounter{ISH}{#2}\addtocounter{ISH}{-1}%
\multido{\i=0+1}{#1}{\psline(\i,0)(\i,\theISH)}%
\multido{\i=0+1}{#2}{\psline(0,\i)(\theISL,\i)}%
}

\newcommand{\colD}[1]{%
\multido{\i=0+1}{#1}{\rput(0,\i){\flecheD}}%
}

\newcommand{\colG}[1]{%
\multido{\i=0+1}{#1}{\rput(1,\i){\flecheG}}%
}

\newcommand{\linH}[1]{%
\multido{\i=0+1}{#1}{\rput(\i,0){\flecheH}}%
}

\newcommand{\linB}[1]{%
\multido{\i=0+1}{#1}{\rput(\i,1){\flecheB}}%
}

\newcommand{\col}[1]{%
\multido{\i=0+1}{#1}{\rput(0,\i){\psline(0,0)(1,0)}}%
}

\newcommand{\halfCC}[5]{%
\degrees[360]%
\multido{\n=#3+#4}{#5}{\rput(#1,#2){\psarc(0,0){\n}{-90}{90}}}%
}

\newcommand{\IceSquare}[1]{%
\rput(1,1){\IceGrid{#1}{#1}}%
\rput(0,1){\colD{#1}}%
\rput(1,0){\linB{#1}}%
\rput(#1,1){\colG{#1}}%
\rput(1,#1){\linH{#1}}%
}

\newcounter{ISdoublesize}
\newcommand{\HTIceEven}[1]{%
\setcounter{ISdoublesize}{#1}\addtocounter{ISdoublesize}{#1}%
\rput(1,1){\IceGrid{#1}{\theISdoublesize}}%
\rput(0,1){\colD{\theISdoublesize}}%
\rput(1,0){\linB{#1}}%
\rput(1,\theISdoublesize){\linH{#1}}%
\rput(#1,#1){\rput(0,.5){\halfCC{0}{0}{.5}{1.0}{#1}}}%
\rput(#1,#1){\psline[linestyle=dotted](.25,.5)(.75,.5)}%
}

\newcommand{\HTIceOdd}[1]{%
\setcounter{ISdoublesize}{#1}\addtocounter{ISdoublesize}{#1}%
\addtocounter{ISdoublesize}{1}%
\rput(1,1){\IceGrid{#1}{\theISdoublesize}}%
\rput(0,1){\colD{\theISdoublesize}}%
\rput(1,0){\linB{#1}}%
\rput(1,\theISdoublesize){\linH{#1}}%
\rput(#1,0){\psline(1,1)(1,#1)\rput(1,1){\flecheB}}%
\psarc(#1,#1){1}{0}{90}%
\rput(#1,#1){\halfCC{1}{1}{1}{1}{#1}}%
\rput(#1,1){\col{#1}}%
\rput(#1,#1){\rput(0,2){\col{#1}}}%
\rput(#1,#1){\psline[linestyle=dotted](.25,.25)(1.25,1.25)}
}

\newcounter{QTsize}

\newcommand{\Hcrossing}{%
\psbezier(0,0)(.4,0)(.6,1)(1,1)%
\psbezier(0,1)(.4,1)(.6,0)(1,0)%
}

\newcommand{\upleft}{%
\begin{pspicture}(.2,.2)\psset{linewidth=.4pt,arrowsize=2.5pt}
\psline(.2,0)(.2,.15)\psarc(.15,.15){.05}{0}{90}\psline{->}(.15,.2)(0,.2)%
\end{pspicture}%
}

\newcommand{\downright}{%
\begin{pspicture}(.2,.2)\psset{linewidth=.4pt,arrowsize=2.5pt}
\psline(0,.2)(.15,.2)\psarc(.15,.15){.05}{0}{90}\psline{->}(.2,.15)(.2,0)%
\end{pspicture}
}

\newcommand{\upc}{%
\begin{pspicture}(.2,.2)\psset{linewidth=.4pt,arrowsize=2.5pt}%
\psline(0,0)(.1,0)\psarc(.1,.1){.1}{270}{90}\psline{->}(.1,.2)(0,.2)%
\end{pspicture}%
}

\newcommand{\downc}{%
\begin{pspicture}(.2,.2)\psset{linewidth=.4pt,arrowsize=2.5pt}%
\psline(0,.2)(.1,.2)\psarc(.1,.1){.1}{270}{90}\psline{->}(.1,0)(0,0)%
\end{pspicture}%
}

\newcommand{\ZhUp}{%
Z_{\textsc{HT}}^{\upc}%
}

\newcommand{\ZhDown}{%
Z_{\textsc{HT}}^{\downc}%
}

\newcommand{\ZhUpleft}{%
Z_{\textsc{HT}}^{\upleft}%
}

\newcommand{\ZhDownright}{%
Z_{\textsc{HT}}^{\downright}%
}

\author{Jean-Christophe Aval}
\thanks{This work has been supported by the ANR project MARS (BLAN06-2$\_$0193)} 

\address{Jean-Christophe Aval, LaBRI, Universit\'e Bordeaux 1, CNRS\\ 351 cours
 de la Lib\'eration, 33405 Talence cedex, FRANCE}

\title{On the symmetry of the partition function of some square ice models}

\date{\today}

\begin{document}

\maketitle

\begin{abstract}
We consider the partition function $Z(N;x_1,\dots,x_N,y_1,\dots,y_N)$ of the square ice model with domain wall boundary. We give a simple proof of the symmetry of $Z$ with respect to all its variables when the global parameter $a$ of the model is set to the special value $a=\exp(i\pi/3)$. Our proof does not use any determinantal interpretation of $Z$ and can be adapted to other situations (for examples to some symmetric ice models).
\end{abstract}

\section{Introduction}\label{sec:intro}

An {\em alternating sign matrix} (ASM) is a square matrix with entries in $\{-1,0,1\}$ and such that in any row and column: the non-zero entries alternate in sign, and their sum is equal to $1$. Their numbers appear in the so-called {\em Razumov-Stroganov conjecture} related to the $O(1)$ loop model \cite{RS1,RS2,cea}. Their enumeration formula was conjectured by Mills, Robbins and Rumsey \cite{MRR}, and proved by Zeilberger \cite{zeil}, and almost simultaneously by Kuperberg \cite{kup1}. Kuperberg used a bijection between the ASM's and the states of a statistical square ice model, for which he studied and computed the partition function. He also used these tools in \cite{kup} to obtain many enumeration or equinumeration results for various classes of symmetries of ASM's, most of them having been conjectured by Robbins \cite{robbins}. The same method was recently used to obtain the enumeration of ASM's invariant \cite{RSQT} or quasi-invariant \cite{AD} under a quarter-turn rotation. A property useful in all these works, which was established by Stroganov \cite{IKdet}, states that the partition function $Z(N;x_1,\dots,x_N,y_1,\dots,y_N)$ of the (unrestrited) square ice model with domain wall boundary is symmetric in all its variables when the global parameter $a$ of the model is set to the special value $a=\exp(i\pi/3)$. We give here another proof of this result. This proof is somehow ``elementary'' since it does not use any determinantal interpretation. Moreover, it can be adapted to other cases, for example to some symmetric models, to models in which there can be lines that carry more than one spectral parameter, or in which some fixed oriented edges are specified.

This paper is organized as follows: in Section \ref{sec:def}, we recall the definitions of the ice models and of their associated weights; in Section \ref{sec:result}, we give and prove the main result, \ie the symmetry of $Z$; in the last section, we present how we can obtain in the same way symmetry properties for other ice models.

\section{Definitions}\label{sec:def}

\subsection{Notations}

We recall here the main definitions and refer to \cite{kup} for details and many examples. 

Let $a\in\mathbb{C}$ be a global parameter.
For any complex number $x$ different from zero, we denote $\overline{x}=1/x$, and we define:
\begin{equation}
  \sigma(x) =  x - \overline{x}.
\end{equation}

If $G$ is a tetravalent graph, an {\em ice state} of $G$ is an orientation of the edges such that every tetravalent vertex has exactly two incoming and two outcoming edges.

A parameter $x\neq 0$ is assigned to any tetravalent vertex of the graph $G$. 
Then this vertex gets a weight, which depends on its orientations, as shown on 
Figure~\ref{fig:poids_6V}.

\begin{figure}[htbp]
  \begin{center}
\psset{unit=0.5cm}
\begin{pspicture}[.5](3,3)
\rput(1.5,1.5){
   \psline(-1,0)(1,0)
   \psline(0,-1)(0,1)
   \psdot[dotsize=.2]
   \rput(-.5,-.5){$x$}
  }
\end{pspicture}
$=$
    \begin{pspicture}[.5](19,5)
\rput(2,2.5){%
    \rput(-1,0){\flecheD}%
    \rput(1,0){\flecheG}%
    \rput(0,0){\flecheH}%
    \rput(0,0){\flecheB}%
    \rput[t](0,-1.2){$\sigma(a^2)$}%
    \rput[t](0,-2.2){$1$}%
}
\rput(5,2.5){%
    \rput(0,0){\flecheG}%
    \rput(0,0){\flecheD}%
    \rput(0,-1){\flecheH}%
    \rput(0,1){\flecheB}%
    \rput[t](0,-1.2){$\sigma(a^2)$}%
    \rput[t](0,-2.2){$-1$}%
}
\rput(8,2.5){%
    \rput(0,0)\flecheG
    \rput(0,0)\flecheH
    \rput(1,0)\flecheG
    \rput(0,-1)\flecheH
    \rput[t](0,-1.2){$\sigma(ax)$}
    \rput[t](0,-2.2){$0$}%
}
\rput(11,2.5){
    \rput(0,0)\flecheD
    \rput(0,0)\flecheB
    \rput(-1,0)\flecheD
    \rput(0,1)\flecheB
    \rput[t](0,-1.2){$\sigma(ax)$}
    \rput[t](0,-2.2){$0$}%
}
\rput(14,2.5){
    \rput(-1,0)\flecheD
    \rput(0,0)\flecheD
    \rput(0,-1)\flecheH
    \rput(0,0)\flecheH
    \rput[t](0,-1.2){$\sigma(a\overline{x})$}
    \rput[t](0,-2.2){$0$}%
}
\rput(17,2.5){
    \rput(0,0)\flecheG
    \rput(1,0)\flecheG
    \rput(0,0)\flecheB
    \rput(0,1)\flecheB
    \rput[t](0,-1.2){$\sigma(a\overline{x})$}
    \rput[t](0,-2.2){$0$}%
}
    \end{pspicture}
  \end{center}
\caption{The 6 possible orientations, their associated weights, and the corresponding entries in ASM's}
\label{fig:poids_6V}
\end{figure}

It is sometimes easier to assign parameters, not to each vertex of the graph, but to the lines that compose the graph. In this case, the weight of a vertex is defined as:

\begin{displaymath}
\psset{unit=.5cm}
  \begin{pspicture}[.45](2,2)
\psline(1,0)(1,2)\psline(0,1)(2,1)
\rput[r](-.2,1){$x$}\rput[t](1,-.2){$y$}
  \end{pspicture}\ =\ 
  \begin{pspicture}[.45](2,2)
\psline(0,1)(2,1)\psline(1,0)(1,2)
\rput(.25,.25){$x\overline{y}$}
  \end{pspicture}
\end{displaymath}

\bigskip

When this convention is used, a parameter explicitly written at a vertex replaces the quotient of the parameters of the lines.

The partition function of a given ice model is then defined as the summation over all its states of the product of the weights of the vertices.

To simplify notations, we will denote by $X_{N}$ the vector of variables
$(x_1,\dots, x_N)$. We use the notation
$X\backslash x$ to denote the vector $X$ without the variable $x$.

\subsection{ASM's and square ice model}

We give in Figure \ref{fig:Z}
the ice model corresponding to (unrestricted) ASM's and its partition functions. The bijection between ASM's and states of the square ice model with ``domain wall boundary'' is now well-known (\cf \cite{kup}). The correspondence between orientations of the ice model and entries of ASM's is given in Figure \ref{fig:poids_6V}.

\begin{figure}[htbp]
  \begin{displaymath}
    Z(N;x_1,\dots, x_{N}, y_{1},\dots, y_{N}) =
    \psset{unit=.5cm}
    \begin{pspicture}[.5](7,7)
      \rput(1,0){\IceSquare{6}}
      \rput[r](1,1){$x_{1}$}\rput[r](1,2){$x_{2}$}\rput[r](1,6){$x_{N}$}
      \rput[t](2,-.1){$y_{1}$}\rput[t](3,-.1){$y_2$}\rput[t](7,-.1){$y_{N}$}
    \end{pspicture}
  \end{displaymath}
\caption{Partition function for ASM's of size $N$}
\label{fig:Z}
\end{figure}

\subsection{Yang-Baxter equation}

To deal with partition functions of ice models, {\bf the} crucial tool is Yang-Baxter equation, that we recall below.
\begin{lemm}\label{lemm:YB}{\em [Yang-Baxter equation]}
  If $xyz=\overline{a}$, then
  \begin{equation}
    \begin{pspicture}[.5](2,2)
      \SpecialCoor
      \rput(1,1){
	\psarc(-1.732,0){2}{330}{30}
        \psarc(.866,1.5){2}{210}{270}
        \psarc(.866,-1.5){2}{90}{150}
	\rput(1;60){$x$}
        \rput(1;180){$y$}
        \rput(1;300){$z$}
      }
    \end{pspicture}
     = 
     \begin{pspicture}[.5](2,2)
      \rput(1,1){
	\SpecialCoor
	\psarc(1.732,0){2}{150}{210}
        \psarc(-.866,1.5){2}{270}{330}
        \psarc(-.866,-1.5){2}{30}{90}
        \rput(1;240){$x$}
        \rput(1;0){$y$}
        \rput(1;120){$z$}
      }
     \end{pspicture}.
\label{eq:YB}
  \end{equation}
\end{lemm}

\section{Main result}\label{sec:result}

The following result has been obtained by Stroganov \cite{IKdet}.
\begin{theo}\label{theo:main}{\em (Stroganov)}
When $a=\omega_6=\exp(i\pi/3)$, the partition function $Z(N;X_N,Y_N)$ is symmetric in {\bf all} its variables $x_1,\dots,x_N,y_1,\dots,y_N$.
\end{theo}

To prove this result Stroganov \cite{IKdet} uses a determinantal interpretation of $Z$. We want here to give a proof that only uses Yang-Baxter equation to study the partition function.

The method used was introduced by Kuperberg \cite{kup}: observe that $Z$ is a Laurent polynomial, then give enough specialization of one of its variable to imply the desired property.

\subsection{Laurent polynomial}

Since the weight of any vertex is a Laurent polynomial in the variables $x_i$'s and $y_i$'s, the partition function $Z$ is a Laurent polynomial in these variables. Moreover it is a centered Laurent polynomial, \ie its lowest degree is the opposite of its highest degree (called the half-width of the polynomial). Since any row and column of an ASM has at least one non-zero entry, which corresponds to a constant $\si(a^2)$, we get the following property.
\begin{lemm}\label{lemm:laurent}
The partition function $Z(N;X_N,Y_N)$ is a Laurent polynomial in any of its variables of half-width $N-1$ and of parity the parity of $N-1$. 
\end{lemm}

\subsection{Partial symmetry}
The following lemma gives a (now classical) example of use of the Yang-baxter equation.

\begin{lemm}\label{lemm:echange_lignes}
  \begin{equation}
    \psset{unit=.5cm}
    \begin{pspicture}[.5](5.5,2)
      \rput[r](1,.5){$x$}\rput[r](1,1.5){$y$}
      \rput(1,.5){\flecheD}\rput(1,1.5){\flecheD}
      \psline(2,0)(2,2)\psline(3,0)(3,2)\psline(4.5,0)(4.5,2)
      \psline(2,.5)(3.5,.5)\psline(2,1.5)(3.5,1.5)
      \psline(4,.5)(4.5,.5)\psline(4,1.5)(4.5,1.5)
      \rput(3.75,1){$\dots$}
      \rput(5.5,.5){\flecheG}
      \rput(5.5,1.5){\flecheG}
    \end{pspicture} = 
    \begin{pspicture}[.5](5.5,2)
      \rput[r](1,.5){$y$}\rput[r](1,1.5){$x$}
      \rput(1,.5){\flecheD}\rput(1,1.5){\flecheD}
      \psline(2,0)(2,2)\psline(3,0)(3,2)\psline(4.5,0)(4.5,2)
      \psline(2,.5)(3.5,.5)\psline(2,1.5)(3.5,1.5)
      \psline(4,.5)(4.5,.5)\psline(4,1.5)(4.5,1.5)
      \rput(3.75,1){$\dots$}
      \rput(5.5,.5){\flecheG}
      \rput(5.5,1.5){\flecheG}      
    \end{pspicture}.
    \label{eq:symetrie_Z}
  \end{equation}
\end{lemm}

\proof
We multiply the left-hand side by $\sigma(a\overline{z})$, with
  $z=\overline{a}x\overline{y}$. We get
  \begin{eqnarray*} \psset{unit=.5cm}
    \sigma(a\overline{z})
    \begin{pspicture}[.5](5.5,2)
      \rput[r](1,.5){$x$}\rput[r](1,1.5){$y$}
      \rput(1,.5){\flecheD}\rput(1,1.5){\flecheD}
      \psline(2,0)(2,2)\psline(3,0)(3,2)\psline(4.5,0)(4.5,2)
      \psline(2,.5)(3.5,.5)\psline(2,1.5)(3.5,1.5)
      \psline(4,.5)(4.5,.5)\psline(4,1.5)(4.5,1.5)
      \rput(3.75,1){$\dots$}
      \rput(5.5,.5){\flecheG}
      \rput(5.5,1.5){\flecheG}
    \end{pspicture} & = & \psset{unit=.5cm}
    \begin{pspicture}[.5](6.5,2)
      \rput[r](1,.5){$y$}\rput[r](1,1.5){$x$}
      \rput(1,.5){\flecheD}\rput(1,1.5){\flecheD}
      \rput(2,.5){\Hcrossing}\rput[r](2.4,1){$z$}
      \rput(1,0){
	\psline(2,0)(2,2)\psline(3,0)(3,2)\psline(4.5,0)(4.5,2)
	\psline(2,.5)(3.5,.5)\psline(2,1.5)(3.5,1.5)
	\psline(4,.5)(4.5,.5)\psline(4,1.5)(4.5,1.5)
	\rput(3.75,1){$\dots$}
	\rput(5.5,.5){\flecheG}
	\rput(5.5,1.5){\flecheG}}
    \end{pspicture} \\ 
      & = & 
    \psset{unit=.5cm}
    \begin{pspicture}[.5](5.5,2)
      \rput[r](1,.5){$y$}\rput[r](1,1.5){$x$}
      \rput(1,.5){\flecheD}\rput(1,1.5){\flecheD}
      \psline(2,0)(2,2)
      \rput(2,.5){\Hcrossing}\rput[r](2.4,1){$z$}
      \psline(3,0)(3,2)\psline(4.5,0)(4.5,2)
      \psline(3,.5)(3.5,.5)\psline(3,1.5)(3.5,1.5)
      \rput(3.75,1){$\dots$}
      \psline(4,.5)(4.5,.5)\psline(4,1.5)(4.5,1.5)
      \rput(5.5,.5){\flecheG}
      \rput(5.5,1.5){\flecheG}
    \end{pspicture} \\
      & = & 
    \psset{unit=.5cm}
    \begin{pspicture}[.5](6.5,2)
      \rput[r](1,.5){$y$}\rput[r](1,1.5){$x$}
      \rput(1,.5){\flecheD}\rput(1,1.5){\flecheD}
      \psline(2,0)(2,2)
      \psline(3,0)(3,2)
      \psline(2,.5)(3.5,.5)\psline(2,1.5)(3.5,1.5)
      \rput(3.75,1){$\dots$}
      \psline(4,.5)(4.5,.5)\psline(4,1.5)(4.5,1.5)
      \rput(4.5,.5){\Hcrossing}\rput[r](4.9,1){$z$}
      \psline(5.5,0)(5.5,2)
      \rput(6.5,.5){\flecheG}
      \rput(6.5,1.5){\flecheG}
    \end{pspicture} \\
      & = & 
    \psset{unit=.5cm}
    \begin{pspicture}[.5](6.5,2)
      \rput[r](1,.5){$y$}\rput[r](1,1.5){$x$}
      \rput(1,.5){\flecheD}\rput(1,1.5){\flecheD}
      \psline(2,0)(2,2)\psline(3,0)(3,2)\psline(4.5,0)(4.5,2)
      \psline(2,.5)(3.5,.5)\psline(2,1.5)(3.5,1.5)
      \psline(4,.5)(4.5,.5)\psline(4,1.5)(4.5,1.5)
      \rput(3.75,1){$\dots$}
      \rput(4.5,.5){\Hcrossing}\rput[l](5.1,1){$z$}
      \rput(6.5,.5){\flecheG}
      \rput(6.5,1.5){\flecheG}      
    \end{pspicture} \\
      & = & 
    \psset{unit=.5cm}
    \begin{pspicture}[.5](5.5,2)
      \rput[r](1,.5){$y$}\rput[r](1,1.5){$x$}
      \rput(1,.5){\flecheD}\rput(1,1.5){\flecheD}
      \psline(2,0)(2,2)\psline(3,0)(3,2)\psline(4.5,0)(4.5,2)
      \psline(2,.5)(3.5,.5)\psline(2,1.5)(3.5,1.5)
      \psline(4,.5)(4.5,.5)\psline(4,1.5)(4.5,1.5)
      \rput(3.75,1){$\dots$}
      \rput(5.5,.5){\flecheG}
      \rput(5.5,1.5){\flecheG}      
    \end{pspicture} \sigma(a\overline{z})
  \end{eqnarray*}
\endproof

As a consequence, we get the ``partial symmetry'' of $Z$, which is true whatever the value of the global parameter $a$.
\begin{lemm}\label{lemm:sym}
The functions $Z(N;X_{N},Y_N)$  is symmetric separately in the two sets of variables $X_N$ and $Y_N$.
\end{lemm}

\subsection{Specialization}
\begin{lemm}{\em [specialization of $Z$; Kuperberg]}
\label{lemm:specZ}
If we denote
\begin{eqnarray*}
  A(y_{1},X_{N}\backslash x_1,Y_N\backslash y_{1}) & = & 
    \prod_{2\leq k\leq N} \sigma(a x_k \overline{y}_{1})
    \prod_{1\leq k\leq N} \sigma(a^2 y_{1} \overline{y}_k),\\
  \overline{A}(y_{1},X_{N}\backslash x_1,Y_N\backslash y_{1}) & = & 
    \prod_{2\leq k\leq N} \sigma(a y_{1} \overline{x}_k) 
    \prod_{1\leq k\leq N} \sigma(a^2 y_k \overline{y}_{1}),
\end{eqnarray*}
then we have:
  \begin{eqnarray}
 Z(N;{\bf \overline{a}y_{1}},X_{N}\backslash x_1,Y_N) & = & 
 \overline{A}(y_{1},X_{N}\backslash x_1,Y_N\backslash y_{1}) Z(N-1;X_{N}\backslash x_1,Y\backslash y_{1}) \label{eq:Zbax}\\
 Z(N;{\bf ay_{1}},X_{N}\backslash x_1,Y_N) & = & 
  A(y_{1},X_{N}\backslash x_1,Y_N\backslash y_{1})
 Z(N-1; X_{N}\backslash x_1,Y_N\backslash y_{1}) \label{eq:Zax}.
  \end{eqnarray}
\end{lemm}

\begin{figure}[htbp]
    \begin{center}\psset{unit=.5cm}
      \begin{pspicture}(-1,-1)(19,8)
	\rput(0,0){\IceSquare{6}}
	\rput[r](0,1){$x_1=\overline{a}y_{1}$}\rput[r](0,6){$x_{N}$}
        \rput[r](0,2){$x_2$}
        \rput[t](1,-.1){$y_{1}$}\rput[t](6,-.1){$y_{N}$}
	\multido{\i=1+1}{5}{\rput(1,\i){\flecheH}}
        \multido{\i=2+1}{5}{\rput(1,\i){\flecheD}\rput(\i,1){\flecheG}\rput(\i,2){\flecheB}}

	\rput(12,0){%
	  \IceSquare{6}
          \rput[r](0,6){$x_1=ay_{1}$}\rput[r](0,5){$x_{N}$}\rput[r](0,1){$x_2$}
          \rput[t](1,-.1){$y_{1}$}\rput[t](6,-.1){$y_{N}$}
          \multido{\i=2+1}{5}{\rput(1,\i){\flecheB}\rput(\i,5){\flecheH}
             \rput(\i,6){\flecheG}}
          \multido{\i=1+1}{5}{\rput(1,\i){\flecheD}}
	}
      \end{pspicture}
      \caption{Fixed edges for (\ref{eq:Zbax}) on the left and
      (\ref{eq:Zax}) on the right}
    \label{fig:recZ}
    \end{center}
  \end{figure}
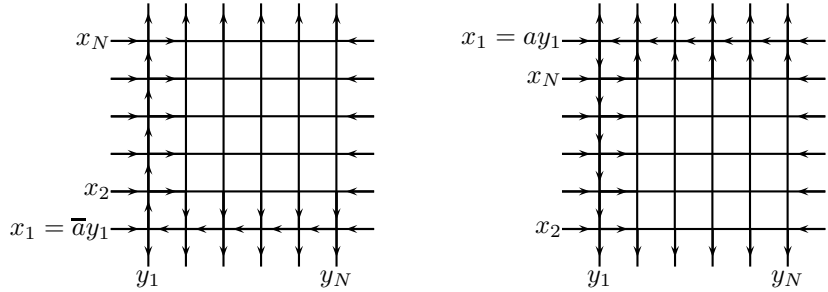

\proof
We recall the method to prove equation \pref{eq:Zbax}. We observe that when $x_1=\bar a y_{1}$, the parameter of the vertex at the crossing of the two lines of parameter $x_1$ and $y_{1}$ is $\bar a$. Thus the weight of this vertex is $\si(a\ba)=\si(1)=0$ unless the orientation of this vertex is the second one on Figure \ref{fig:poids_6V}. But this orientation implies the orientation of all vertices in the row $x_1$ and in the column $y_{1}$, as shown on Figure \ref{fig:recZ}. 
The fixed part gives the partition function $Z$ in size $N-1$, without parameters $x_1$ and
  $y_{1}$, and the weights of the fixed part gives the factor $\overline{A}(\dots)$.

The case of \pref{eq:Zax} is similar, after using Lemma \ref{lemm:sym} to put the line $x_{1}$ at the top of the grid.

\endproof

\subsection{Conclusion}

We are now in a position to complete the proof of Theorem \ref{theo:main}.
From now on we set the global parameter to thz special value $a=\exp(i\pi/3)$. This special value of the parameter implies the following equalities: 
\begin{equation}\label{spec}
\si(a)=\si(a^2) \ \ \ \ \ \ \ \si(a^2x)=-\si(\ba x)=\si(a\bx).
\end{equation}

The proof of Theorem \ref{theo:main} is done by induction on $N$. We easily check the property for $N=1$. Now if $N\ge 2$, because of Lemma \ref{lemm:sym}, we want to obtain the symmetry of $Z$ with respect to $x_2,y_2$.

We use Lemma \ref{lemm:laurent} to reduce the proof of Theorem \ref{theo:main} to the proof of the symmetry of $Z$ in at least $N$ (independent) specializations of the variable $x_1$.
When $x_1=ay_1$,  Lemma \ref{lemm:specZ} gives:
\begin{equation}
 Z(N;{\bf ay_{1}},X_{N}\backslash x_1,Y_N) =  
  A(y_{1},X_{N}\backslash x_1,Y_N\backslash y_{1})
 Z(N-1; X_{N}\backslash x_1,Y_N\backslash y_{1}) \label{eq:Zax2}.
\end{equation}
By recurrence, we have that $ Z(N-1; X_{N}\backslash x_1,Y_N\backslash y_{1})$ is symmetric in $x_2,y_2$. The terms of $A(y_{1},X_{N}\backslash x_1,Y_N\backslash y_{1})$ involving the variables $x_2,y_2$ are:
$$\si(ax_2\bar y_1)\,\si(a^2y_1\bar y_2).$$
We use equation \pref{spec} to write:
$$\si(ax_2\bar y_1)\,\si(a^2y_1\bar y_2)=\si(ax_2\bar y_1)\,\si(ay_2\bar y_1)$$
which is clearly symmetric in the variables $x_2,y_2$. We obtain in the same way the symmetry when $x_1=\bar ay_1$; and by Lemma \ref{lemm:sym} we get this symmetry for the $2(N-1)$ special values $x_1=a^{\pm 1} y_k$ for $k=1,3,\dots,N$, which is more than enough to conclude the proof of Theorem \ref{theo:main}.

\section{Other ice models}\label{sec:ZHT}

The method used to prove the global symmetry of $Z$ may be adapted to other ice models. We illustrate this with half-turn symmetric ASM's (HTASM's). This example shows how our method can be used to prove global symmetries, or partial symmetries (for example when a line of the ice model carries two different spectral parameters -which breaks the homogeneity of the partition function-, or when the orientation of an edge is fixed).

\subsection{HTASM's -- notations and results}

\begin{figure}[htbp]
  \begin{eqnarray*}
    Z_{\textsc{HT}}(2N; x_1,\dots,x_{N-1}(,x,y),y_{1},\dots,y_{N})&  =& 
\psset{unit=.4cm}
    \begin{pspicture}[.5](9,9)
      \rput(1,0){\HTIceEven{4}}
      \rput[r](1,1){$x_{1}$}\rput[r](1,3){$x_{N-1}$}\rput[r](1,5){$y$}
      \rput[r](1,4){$x$}
      \rput[t](2,-.1){$y_{1}$}\rput[t](3,-.1){}\rput[t](5,-.1){$y_{N}$}
    \end{pspicture}\\ 
\psset{unit=.4cm}
\begin{pspicture}[.5](10,10)
  \rput(1,0){\HTIceOdd{4}}
  \rput[r](1,1){$x_{1}$}\rput[r](1,2){$x_{2}$}\rput[r](1,4){$x_{N}$}
  \rput[r](1,5){$x$}
  \rput[t](2,-.1){$y_{1}$}\rput[t](3,-.1){}\rput[t](5,-.1){$y_{N}$}
  \rput[t](6,-.1){$y$}
\end{pspicture} & = & 
Z_{\textsc{HT}}(2N+1;x_1,\dots,x_{N},y_{1},\dots,y_{N},(x,y)) 
  \end{eqnarray*}
\caption{Partition functions for HTASM's}
\label{fig:ZHT}
\end{figure}
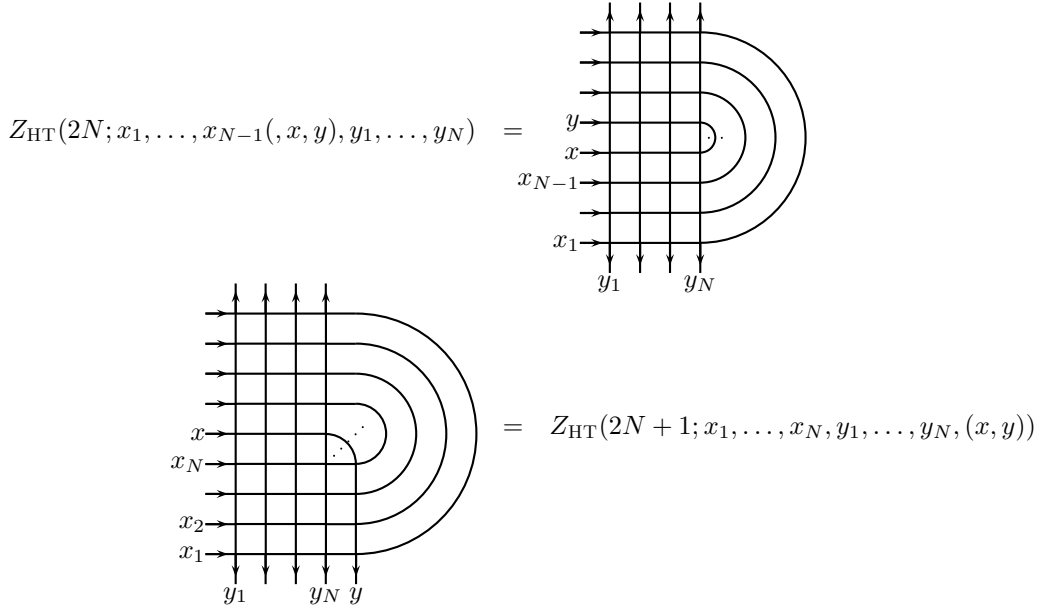

The ice models corresponding to HTASM's were introduced by Kuperberg \cite{kup} for the even size and by Razumov and Stroganov \cite{RSQT} for the odd size. We recall these models on Figure \ref{fig:ZHT}. The spectral parameters on this figure are slightly different from the original ones to better suit the proof. The dotted lines mean a change of parameter: on one side the parameter is $x$, whereas on the other side it is $y$. 

We will prove the following theorem, for the special value $a=\exp(i\pi/3)$, and for the specialization $x=y$ in the even case (which corresponds to the original definition of $Z_{\textsc{HT}}(2N)$.

\begin{theo}\label{theo:mainHT}
When the global parameter $a$ is set to the special value $=\exp(i\pi/3)$, the function $Z_{\textsc{HT}}(2N;X_{N-1},(x_N,x_N),Y_N)$ is symmetric with respect to all its $2N$ variables, and the function $Z_{\textsc{HT}}(2N+1;X_{N},(x,y),Y_N)$ is symmetric with respect to its $2N$ variables $x_1,\dots,x_N,y_1,\dots,y_N$.
\end{theo}

The property on $Z_{\textsc{HT}}(2N+1)$ may be deduced from the main result of \cite{RSQT}. The assertion about $Z_{\textsc{HT}}(2N)$ is new.

\subsection{Proofs}

Since the method is essentially the same as in the case of unrestricted ASM's, we shall give the main steps of the proof and only insist on the difference with the previous case.

\subsubsection{Laurent polynomials}

To deal with Laurent polynomials of given parity in the variable $y$, and thus divide by two the number of required specializations of this variable, 
we shall group together the states with a given orientation (indicated as subscripts in the following notations) at the edge where the parameters $x$ and $y$ meet.

So let us consider the partition functions 
$\ZhUp(2N;X_{N-1},(x,y),Y_N)$ and
$\ZhDown(2N;X_{N-1},(x,y),Y_N)$, respectively parts with the parity of $N-1$ and of $N$ of
$Z_{\textsc{HT}}(2N;X_{N-1},(x,y),Y_N)$ in $y$; and
$\ZhUpleft(2N+1;X_{N},Y_N,(x,y))$ and
$\ZhDownright(2N+1;X_{N},Y_N,(x,y))$, respectively parts with the parity of $N$ and of $N-1$ of
$Z_{\textsc{HT}}(2N+1;X_{N},Y_N,(x,y))$ in $y$.

\begin{lemm}\label{lemm:LaurentHT}
The functions 
$\ZhUp(2N;X_{N-1},(x,y),Y_N)$, 
$\ZhDown(2N;X_{N-1},(x,y),Y_N)$,
$\ZhUpleft(2N+1;X_{N},Y_N,(x,y))$ and
$\ZhDownright(2N+1;X_{N},Y_N,(x,y))$
are centered Laurent polynomials
in the variable $y$,
odd or even, of respective half-widths $N-1$, $N$, $N-1$, and
$N$.
\end{lemm}

\subsubsection{Partial symmetries}

Since Lemma \ref{lemm:echange_lignes} may be easily adapted to HTASM's, we get:

\begin{lemm}\label{lemm:symHT}
The functions 
$\ZhUp(2N;X_{N-1},(x,y),Y_N)$, 
$\ZhDown(2N;X_{N-1},(x,y),Y_N)$,
$\ZhUpleft(2N+1;X_{N},Y_N,(x,y))$ and
$\ZhDownright(2N+1;X_{N},Y_N,(x,y))$
are symmetric separately in the variables $x_i$ and in the variable $y_i$. 

Moreover if we specialize $x=y=x_N$ in the even case, we get the symmetry with respect to the set $X_N$.
\end{lemm}

We now have to deal with the symmetry in the couple $(x,y)$, in the even case.

The easy transformation
  \begin{equation}
    \psset{unit=.5cm}
    \begin{pspicture}[.45](2.5,1)
      \Hcrossing
      \psarc(1,.5){.5}{270}{90}
      \rput[r](.3,.5){$z$}
    \end{pspicture} = \left(\sigma(az)+\sigma(a^2)\right)
    \left(
    \begin{pspicture}[.45](1,1)
      \rput(0,1){\flecheD}\rput(1,0){\flecheG}
    \end{pspicture}\ +\ 
    \begin{pspicture}[.45](1,1)
      \rput(0,0){\flecheD}\rput(1,1){\flecheG}
    \end{pspicture}\right)
    \label{eq:boucle}
  \end{equation}

together with Yang-Baxter equation \pref{eq:YB} gives the following lemma.

\begin{lemm}\label{lemm:echange_boucle}
  \begin{eqnarray}
    \psset{unit=.5cm}
    \begin{pspicture}[.5](5.7,2)
      \rput[r](1,.5){$x$}\rput[r](1,1.5){$y$}
      \rput(1,.5){\flecheD}\rput(1,1.5){\flecheD}
      \psline(2,0)(2,2)\psline(3,0)(3,2)\psline(4.5,0)(4.5,2)
      \psline(2,.5)(3.5,.5)\psline(2,1.5)(3.5,1.5)
      \psline(4,.5)(4.5,.5)\psline(4,1.5)(4.5,1.5)
      \psarc(4.5,1){.5}{270}{90}
      \psline[linestyle=dotted](4.7,1)(5.7,1)
      \rput(3.75,1){$\dots$}
    \end{pspicture} &=&
    \psset{unit=.5cm}
    \frac{\sigma(a^2)+\sigma(x\overline{y})}{\sigma(a^2y\overline{x})}
    \begin{pspicture}[.5](5.7,2)
      \rput[r](1,.5){$y$}\rput[r](1,1.5){$x$}
      \rput(1,.5){\flecheD}\rput(1,1.5){\flecheD}
      \psline(2,0)(2,2)\psline(3,0)(3,2)\psline(4.5,0)(4.5,2)
      \psline(2,.5)(3.5,.5)\psline(2,1.5)(3.5,1.5)
      \psline(4,.5)(4.5,.5)\psline(4,1.5)(4.5,1.5)
      \psarc(4.5,1){.5}{270}{90}
      \psline[linestyle=dotted](4.7,1)(5.7,1)
      \rput(3.75,1){$\dots$}
    \end{pspicture} \label{eq:echange_boucle}\\
    \psset{unit=.5cm}
    \begin{pspicture}[.5](5.7,2)
      \rput[r](1,.5){$x$}\rput[r](1,1.5){$y$}
      \rput(1,.5){\flecheD}\rput(1,1.5){\flecheD}
      \psline(2,0)(2,2)\psline(3,0)(3,2)\psline(4.5,0)(4.5,2)
      \psline(2,.5)(3.5,.5)\psline(2,1.5)(3.5,1.5)
      \psline(4,.5)(4.5,.5)\psline(4,1.5)(4.5,1.5)
      \rput(4.5,1.5){\flecheD}
      \rput(5.5,.5){\flecheG}
      \rput(3.75,1){$\dots$}
    \end{pspicture} &=&
    \psset{unit=.5cm}
\frac{\sigma(x\overline{y})}{\sigma(a^2y\overline{x})}
    \begin{pspicture}[.5](5.7,2)
      \rput[r](1,.5){$y$}\rput[r](1,1.5){$x$}
      \rput(1,.5){\flecheD}\rput(1,1.5){\flecheD}
      \psline(2,0)(2,2)\psline(3,0)(3,2)\psline(4.5,0)(4.5,2)
      \psline(2,.5)(3.5,.5)\psline(2,1.5)(3.5,1.5)
      \psline(4,.5)(4.5,.5)\psline(4,1.5)(4.5,1.5)
      \rput(5.5,1.5){\flecheG}
      \rput(4.5,.5){\flecheD}
      \rput(3.75,1){$\dots$}
    \end{pspicture}  + 
\frac{\sigma(a^2)}{\sigma(a^2y\overline{x})}
    \begin{pspicture}[.5](5.7,2)
      \rput[r](1,.5){$y$}\rput[r](1,1.5){$x$}
      \rput(1,.5){\flecheD}\rput(1,1.5){\flecheD}
      \psline(2,0)(2,2)\psline(3,0)(3,2)\psline(4.5,0)(4.5,2)
      \psline(2,.5)(3.5,.5)\psline(2,1.5)(3.5,1.5)
      \psline(4,.5)(4.5,.5)\psline(4,1.5)(4.5,1.5)
      \rput(4.5,1.5){\flecheD}
      \rput(5.5,.5){\flecheG}
      \rput(3.75,1){$\dots$}
    \end{pspicture} \label{eq:echange_boucle_a}\\
    \psset{unit=.5cm}
    \begin{pspicture}[.5](5.7,2)
      \rput[r](1,.5){$x$}\rput[r](1,1.5){$y$}
      \rput(1,.5){\flecheD}\rput(1,1.5){\flecheD}
      \psline(2,0)(2,2)\psline(3,0)(3,2)\psline(4.5,0)(4.5,2)
      \psline(2,.5)(3.5,.5)\psline(2,1.5)(3.5,1.5)
      \psline(4,.5)(4.5,.5)\psline(4,1.5)(4.5,1.5)
      \rput(4.5,.5){\flecheD}
      \rput(5.5,1.5){\flecheG}
      \rput(3.75,1){$\dots$}
    \end{pspicture} &=&
    \psset{unit=.5cm}
\frac{\sigma(x\overline{y})}{\sigma(a^2y\overline{x})}
    \begin{pspicture}[.5](5.7,2)
      \rput[r](1,.5){$y$}\rput[r](1,1.5){$x$}
      \rput(1,.5){\flecheD}\rput(1,1.5){\flecheD}
      \psline(2,0)(2,2)\psline(3,0)(3,2)\psline(4.5,0)(4.5,2)
      \psline(2,.5)(3.5,.5)\psline(2,1.5)(3.5,1.5)
      \psline(4,.5)(4.5,.5)\psline(4,1.5)(4.5,1.5)
      \rput(5.5,.5){\flecheG}
      \rput(4.5,1.5){\flecheD}
      \rput(3.75,1){$\dots$}
    \end{pspicture}  + 
\frac{\sigma(a^2)}{\sigma(a^2y\overline{x})}
    \begin{pspicture}[.5](5.7,2)
      \rput[r](1,.5){$y$}\rput[r](1,1.5){$x$}
      \rput(1,.5){\flecheD}\rput(1,1.5){\flecheD}
      \psline(2,0)(2,2)\psline(3,0)(3,2)\psline(4.5,0)(4.5,2)
      \psline(2,.5)(3.5,.5)\psline(2,1.5)(3.5,1.5)
      \psline(4,.5)(4.5,.5)\psline(4,1.5)(4.5,1.5)
      \rput(4.5,.5){\flecheD}
      \rput(5.5,1.5){\flecheG}
      \rput(3.75,1){$\dots$}
    \end{pspicture} \label{eq:echange_boucle_b}
  \end{eqnarray}
\end{lemm}

We deduce from this lemma the following property of pseudo-symmetry in $(x,y)$ for the functions $\ZhUp(2N;X_{N-1},(x,y),Y_N)$ and $\ZhDown(2N;X_{N-1},(x,y),Y_N)$.

\begin{lemm}\label{lemm:pseudo_sym}
For $\star=\upc,\downc$ and
  $\square=\downc,\upc$ respectively, we have
\begin{equation}\label{eq:pssym}
    Z_{\textsc{HT}}^\star(2N;X_{N-1},(x,y),Y_N)  = 
    \frac{1}{\sigma(a^2y\overline{x})}
    \sigma(a^2) Z_{\textsc{HT}}^\star(2N;X_{N-1},(y,x),Y_N)
\end{equation}
$$\ \ \ \ \ \ \ \ \ \ \ \ \ \ \ \ \ \ \ \ \ \ \ \    +\,\sigma(x\overline{y}) Z_{\textsc{HT}}^\square(2N;X_{N-1},(y,x),Y_N).$$
\end{lemm}

\subsubsection{Specializations}

We now give specializations of the functions $Z_{\textsc{HT}}$ in the variable $x$ or $y$.

\begin{lemm}{\em [specialization of $Z_{\textsc{HT}}$]}
\label{lemm:specZHT}
If we denote
 \begin{eqnarray*}
   A_{H}^{1}(x_1,X_{N}\backslash x_1,Y_N) & = &    
     \prod_{1\leq k\leq N} \sigma(a^2 x_1 \overline{x}_k)
     \prod_{1\leq k\leq N} \sigma(a y_k\overline{x}_1)\\
   \overline{A}_{H}^{1}(x_1,X_{N}\backslash x_1,Y_N) & = & 
     \prod_{1\leq k\leq N} \sigma(a^2x_k\overline{x}_1)
     \prod_{1\leq k\leq N} \sigma(ax_1\overline{y}_k)\\
   A_{H}^{0}(y_{1},X_{N-1},Y_N\backslash y_{1}) & = & 
     \prod_{1\leq k\leq N-1}\sigma(a x_k \overline{y}_1)
     \prod_{1\leq k\leq N} \sigma(a^2 y_1 \overline{y}_k)\\
   \overline{A}_{H}^{0}(y_{1},X_{N-1},Y_N\backslash y_{1}) & = & 
     \prod_{1\leq k\leq N-1} \sigma(ay_1\overline{x}_k)
     \prod_{1\leq k\leq N} \sigma(a^2y_k\overline{y}_1),
 \end{eqnarray*}
  then for $\star=\downright,\upleft$ and
  $\square=\downc,\upc$ respectively, we have
  \begin{align}
Z_{\textsc{HT}}^{\star}(2N+1;X_{N},Y_N,(x,\mathbf{ax_1}))=&
    A_{H}^{1}(x_1,X_{N}\backslash x_{1},Y_N) 
    Z_{\textsc{HT}}^{\square}(2N;X_{N}\backslash x_1,(x_1,x),Y_N)\label{equa:Zht_impair_ax}\\
Z_{\textsc{HT}}^{\square}(2N+1;X_{N},Y_N,(x,{\bf \overline{a}x_1}))=&
    \overline{A}_{H}^{1}(x_1,X_{N}\backslash x_1,Y_N)
     Z_{\textsc{HT}}^{\star}(2N;X_{N}\backslash x_1,(x,x_1),Y_N)\label{equa:Zht_impair_bax}\\
Z_{\textsc{HT}}^{\star}(2N;X_{N-1},(x,{\bf ay_1}),Y_N) =&
  \sigma(ax\overline{y}_{1}) A_{H}^{0}(y_{1},X_{N-1},Y_N\backslash y_{1})
  Z_{\textsc{HT}}^{\square}(2N-1;X_{N-1},Y_N\backslash y_1,(x,y_1)) \label{equa:Zht_pair_ax}\\
Z_{\textsc{HT}}^{\square}(2N;X_{N-1},({\bf \overline{a}y_1},y),Y_N) =&
  \sigma(ay_1\overline{y}) 
  \overline{A}_{H}^{0}(y_{1},X_{N-1},Y_N\backslash y_{1})
   Z_{\textsc{HT}}^{\star}(2N-1;X_{N-1},Y_N\backslash y_1,(y,y_1))\label{equa:Zht_pair_bax}
  \end{align}
\end{lemm}

\proof
The method is almost the same as the one used to prove Lemma \ref{lemm:specZ}. An extended proof may be found in \cite{AD}.
\endproof

\subsubsection{Conclusion}

We are now in a position to conclude the proof of Theorem \ref{theo:mainHT}. From now on we set the global parameter $a$ to the special value $a=\exp(i\pi/3)$. 

Lemma \ref{lemm:symHT} allows us to reduce the proof of the theorem to the proof of the following assertion: the functions $\ZhUp(2N;X_{N-1},(x,y),Y_N)$, 
$\ZhDown(2N;X_{N-1},(x,y),Y_N)$,
$\ZhUpleft(2N+1;X_{N},Y_N,(x,y))$ and
$\ZhDownright(2N+1;X_{N},Y_N,(x,y))$
are symmetric with respect to $x_2,y_2$.

The proof is done by induction on $N$. We easily check the property for small values of $N$. 
We suppose the property true in size $2N$, and we consider the odd case $2N+1$.

Because of Lemma \ref{lemm:LaurentHT}, we have to check the property for enough specializations of the variable $y$. We use Lemma \ref{lemm:specZHT}: equations \pref{equa:Zht_impair_ax} and \pref{equa:Zht_impair_ax} give us the specializations $y=a^{\pm 1} x_1$. 
As in the case of $Z$, the recurrence gives us the symmetry of the factor $Z_{\textsc{HT}}(2N)$, and the symmetry of the factor $A$ is settled through \pref{spec} which implies:
$$\si(ax_2\bar y_1)\,\si(a^2y_1\bar y_2)=\si(ax_2\bar y_1)\,\si(ay_2\bar y_1)$$
the right-hand being clearly symmetric in $x_2,y_2$.

Using Lemma \ref{lemm:symHT} we get $2(N-1)$ specializations satisfying the symmetry, which is enough to imply the symmetry in full generality as soon as $N>2$.

Now we use the property in size $2N+1$ to prove it in size $2N+2$. The method is the same, with the only difference that we have to use Lemma \ref{lemm:echange_boucle} because equation \pref{equa:Zht_pair_bax} gives a specialization for the variable $x$.

The proof of Theorem \ref{theo:mainHT} is now complete.

\bigskip
{\bf Aknowledgment.} The author would like to thank P. Duchon for valuable discussions and for the use of his figures of ice models.


\end{document}